\documentclass[a4paper,12pt]{article}
\usepackage{amsmath}
\usepackage{amsfonts}
\usepackage{amssymb}
\usepackage{latexsym}
\usepackage{epsfig}
\usepackage{graphicx}
\usepackage{oldgerm}
\usepackage{theorem}

\def\N{\mathbb{N}}
\def\C{\mathbb{C}}
\def\R{\mathbb{R}}
\def\W{\mathbb{W}}
\def\Z{\mathbb{Z}}

\def\P{\mathbb{P}}
\def\mbK{\mathbb{K}}
\def\mbA{\mathbb{A}}
\def\supp{{\rm supp}\ }

\def\x{\mib{x}}
\def\y{\mib{y}}

\def\rA{{\rm A}}
\def\rC{{\rm C}}

\def\1{{\bf 1}}

\def\bE{{\bf E}}
\def\bP{{\bf P}}

\def\X{\mib{X}}

\def\Pf{{\rm Pf}}

\theorembodyfont{\itshape}
\newtheorem{thm}{Theorem}[section]
\newtheorem{lem}[thm]{Lemma}

\newcommand{\mib}[1]{\mbox{\boldmath $#1$}}
\newcommand{\SSC}[1]{\section{#1}\setcounter{equation}{0}}
\newcommand{\qed}{\hbox{\rule[-2pt]{3pt}{6pt}}}

\begin{document}

\title{\bf 
Determinantal process starting from
an orthogonal symmetry is a Pfaffian process}
\author{
Makoto Katori
\footnote{
Department of Physics,
Faculty of Science and Engineering,
Chuo University, 
Kasuga, Bunkyo-ku, Tokyo 112-8551, Japan;
e-mail: katori@phys.chuo-u.ac.jp
}}
%%%%%%%%%%%%%%%%%%%%%%%%%%%%%%%%
\date{3 October 2011}
%%%%%%%%%%%%%%%%%%%%%%%%%%%%%%%%
\pagestyle{plain}
\maketitle
\begin{abstract}

When the number of particles $N$ is finite,
the noncolliding Brownian motion (BM) and 
the noncolliding squared Bessel process with
index $\nu > -1$ (BESQ$^{(\nu)}$) are
determinantal processes for arbitrary fixed initial
configurations.
In the present paper we prove that,
if initial configurations are distributed with
orthogonal symmetry,
they are Pfaffian processes in the sense that
any multitime correlation functions are
expressed by Pfaffians.
The $2 \times 2$ skew-symmetric matrix-valued
correlation kernels of the Pfaffians processes
are explicitly obtained by the equivalence
between the noncolliding BM and
an appropriate dilatation of a time reversal
of the temporally inhomogeneous version of
noncolliding BM with finite duration
in which all particles start from the origin, $N \delta_0$,
and by the equivalence between the noncolliding BESQ$^{(\nu)}$
and that of the noncolliding squared generalized meander
starting from $N \delta_0$.

\end{abstract}

%%%%%%%%%%%%%%%%%%%%%%%%%%%%%%
%\newpage
%%%%%%%%%%%%%%%%%%%%%%%
%\footnotesize 
%\tableofcontents
%%%%%%%%%%%%%%%%%%%%%%%
%\vspace{3mm}
\normalsize
%\setlength{\baselineskip}{1cm}

%%%%%%%%%%%%%%%%%%%%%%%%%%%%%%%%%%%%%%%%%%%%%%%%%%%%%%%%%%
%%%  SEC1 %%%%%%%%%%%%%%%%%%%%%%%%%%%%%%%%%%%%%%%%%%%
%%%%%%%%%%%%%%%%%%%%%%%%%%%%%%%%%%%%%%%%%%%%%%%%%%%%%%%%%%
\SSC{Introduction}%%%%%%%%
%%%%%%%%%%%%%%%%%%%%%%%%%%%%%%%%%%%%%%%%%%%%%%%%%%%%%%%%%%

We consider one-dimensional particle systems called
{\it the noncolliding Brownian motion} and
{\it the noncolliding squared Bessel process}.
The former is equivalent with Dyson's Brownian motion (BM) model
with parameter $\beta=2$,
which was introduced as the eigenvalue process
of the Hermitian-matrix-valued BM \cite{Dys62}
corresponding to the Gaussian unitary ensemble (GUE) 
of random matrices \cite{Meh04,For10}.
The latter is a one-parameter family indexed by $\nu > -1$
and abbreviated as noncolliding BESQ$^{(\nu)}$
\cite{KT11}.
When the number of particles is finite,
$N \in \N \equiv \{1,2,3, \dots\}$ and
$\nu \in \N_0 \equiv \N \cup \{0\}$,
the noncolliding BESQ$^{(\nu)}$ 
realizes the eigenvalue process of the 
matrix-valued diffusion process
called the Laguerre process (or complex Wishart process)
\cite{KO01},
whose distribution at each time describes 
squares of singular values of 
$(N+\nu) \times N$ random matrices
in the chiral Gaussian unitary ensemble (chGUE) 
\cite{VZ93,Ver94}, 
and when $\nu=1/2$ (resp. $\nu=-1/2$),
it expresses the stochastic evolution \cite{KT04} of
the squares of positive eigenvalues of
$2N \times 2N$ random matrices in the Gaussian ensemble
of class C (resp. class D) \cite{AZ96,AZ97}.
(See \cite{FNH99,BHL99,KTNK03,Nag03,KT07a,
TW07,SMCR08,KMW09,BFPSW09,RS10,Ols10,
DKZ10,KMW10,RS11} 
for related interacting particle systems.)

In the previous papers \cite{KT09,KT10,KT11}, it was proved that
if the number of particles is finite, 
$N \in \N$,
these two processes are {\it determinantal}
for arbitrary initial configurations
$\xi(\cdot)=\sum_{j=1}^{N}\delta_{x_j}(\cdot)$ with 
$x_1 \leq x_2 \leq \cdots \leq x_N, x_j \in \Lambda, 1 \leq j \leq N$,
where $\delta_y(\cdot)$ denotes 
the delta measure on $y$;
$\delta_y(x)=\delta_{x y}$, 
$\Lambda=\R$ for the noncolliding BM and
$\Lambda=\R_+ \equiv \{x \in \R: x \geq 0\}$
for the noncolliding BESQ$^{(\nu)}$.
Here, given a fixed initial configuration
$\xi$, a process is said to be determinantal,
if there is a function $\mbK(s,x;t,y)$
such that it is continuous with respect to
$(x,y)$ for any fixed $(s, t) \in [0, \infty)^2$
and any multitime correlation function
is given by a determinant in the form
\begin{equation}
\rho^{\xi}
 (t_1,\x^{(1)}_{N_1}; \dots;t_M,\x^{(M)}_{N_M} ) 
=\det_{
\substack{1 \leq j \leq N_{m}, 1 \leq k \leq N_{n} \\ 1 \leq m, n \leq M}
}
[\mbK(t_m, x_{j}^{(m)}; t_n, x_{k}^{(n)} )],
\label{eqn:det1}
\end{equation}
$M \in \N,
0 < t_1 < \cdots < t_M < \infty,
1 \leq N_m \leq N, 1 \leq m \leq M$,
where 
$\x^{(m)}_{N_m}=(x^{(m)}_1, \dots, x^{(m)}_{N_m})$
denotes the points at which observation is performed
at time $t_m$, $1 \leq m \leq M$ \cite{KT07b}.
The function $\mbK$ is called 
the {\it correlation kernel} and it determines 
finite dimensional distributions of the process
through (\ref{eqn:det1}).
For a configuration $\xi(\cdot)=\sum_{j=1}^{N} \delta_{x_j}(\cdot)$,
shift by $w \in \C$ is denoted by
$\tau_{w} \xi(\cdot)=\sum_{j=1}^{N} \delta_{x_j+w}(\cdot)$
and dilatation by factor $c>0$ is denoted by
$c \circ \xi(\cdot)=\sum_{j=1}^{N} \delta_{c x_j}(\cdot)$.
The correlation kernels for the noncolliding BM
and the noncolliding BESQ$^{(\nu)}$ are 
respectively given by
\begin{eqnarray}
\mbK^{\xi}(s,x;t,y) &=& \lim_{\varepsilon \downarrow 0}
\frac{1}{2 \pi i} \int_{\R \setminus[-\varepsilon, \varepsilon]} du
\oint_{\rC_{i u}(\xi)} dz \,
p(s,x|z) 
\frac{\Pi_{\tau_{-z}\xi}(i u-z)}{i u-z} p(-t,i u|y) 
\nonumber\\
\label{eqn:K1}
&&
- \1(s>t) p(s-t,x|y),
\quad (x,y) \in \R^2, (s,t) \in [0, \infty)^2,
\\
\mbK^{\xi}_{\nu}(s,x;t,y) &=& \lim_{\varepsilon \downarrow 0}
\frac{1}{2 \pi i} \int_{-\infty}^{-\varepsilon} du
\oint_{\rC_{u}(\xi)} dz \,
p^{(\nu)}(s,x|z)
\frac{\Pi_{\tau_{-z}\xi}(u-z)}{u-z} p^{(\nu)}(-t, u|y)
\nonumber\\
&&
- \1(s>t) p^{(\nu)}(s-t,x|y),
\quad (x,y) \in (0,\infty)^2, (s,t) \in [0, \infty)^2,
\nu > -1,
\nonumber\\
\label{eqn:K2}
\end{eqnarray}
where $i=\sqrt{-1}$,
$\rC_{z'}(\xi)$ denotes a closed contour on the
complex plane $\C$ encircling the points in 
$\supp \xi \equiv \{x \in \Lambda: \xi(\{x\}) > 0\}$ 
once in the positive direction
but not the point $z'$,
$p$ and $p^{(\nu)}$ are the extended versions
of transition probability densities
of the one-dimensional standard BM \cite{KT10}
and the BESQ$^{(\nu)}$ \cite{KT11},
\begin{equation}
p(t, y|x)
= \left\{ \begin{array}{ll}
\displaystyle{
\frac{1}{\sqrt{2 \pi |t|}} 
e^{-(x-y)^2/2t}},
& \quad t \in \R \setminus \{0\}, x,y \in \C, \cr
\delta(y-x),
& \quad t=0, x,y \in \C,
\end{array} \right.
\label{eqn:p1}
\end{equation}
\begin{equation}
p^{(\nu)}(t, y|x)
=  \left\{ \begin{array}{ll}
\displaystyle{
\frac{1}{2|t|} \left( \frac{y}{x} \right)^{\nu/2}
e^{-(x+y)/2t}
I_{\nu} \left( \frac{\sqrt{xy}}{t} \right)},
&t \in \R \setminus \{0\}, x \in \C \setminus \{0\}, y \in \C, \cr
\displaystyle{
\frac{y^{\nu}e^{-y/2t} }{(2 |t|)^{\nu+1} \Gamma(\nu+1)} },
& t \in \R \setminus \{0\}, x=0, y \in \C, \cr
& \cr
\delta(y-x),
& t=0, x,y \in \C,
\end{array} \right.
\label{eqn:pnu1}
\end{equation}
with the Gamma function 
$\Gamma(z)=\int_0^{\infty} e^{-u} u^{z-1} du$
and the modified Bessel function
$I_{\nu}(z)=\sum_{n=0}^{\infty}(z/2)^{2n+\nu}/
\{\Gamma(n+1) \Gamma(n+1+\nu) \}$, 
$\Pi_{\xi}$ is an entire function
having $\supp \xi$ as the zero set 
expressed by the following Weierstrass canonical product
with genus 0 \cite{Lev96,Nog98},
\begin{equation}
\Pi_{\xi}(z) =\prod_{x \in \xi}
\left(1-\frac{z}{x} \right)
\equiv \prod_{x \in \supp\xi}
\left(1-\frac{z}{x} \right)^{\xi(\{x\})},
\quad z \in \C,
\label{eqn:Phi}
\end{equation}
and $\1(\omega)$ is the indicator of a condition $\omega$;
$\1(\omega)=1$ if $\omega$ is satisfied
and $\1(\omega)=0$ otherwise.
In (\ref{eqn:pnu1})
we have defined $z^{\nu}$ 
to be $\exp(\nu \log z)$, where the argument of $z$
is given its principal value;
$
z^{\nu}=\exp[ \nu \{
\log |z| + \sqrt{-1} {\rm arg} (z)\}],
-\pi < {\rm arg} (z) \leq \pi.
$
We say that the correlation kernels, which are asymmetric,
$\mbK(s,x;t,y) \not= \mbK(t,y, s,x)$
for $s \not= t$, as
(\ref{eqn:K1}) and (\ref{eqn:K2}), are
{\it of Eynard-Mehta type}
\cite{EM98,NF98,BR04,KT07b,KT10b}.

For $N \in \N, \xi=\sum_{j=1}^N \delta_{x_j},
\x=(x_1, \dots, x_N), \beta \geq 1, a > -1, \sigma^2 >0$, let
\begin{eqnarray}
\label{eqn:mub}
&& \mu_{N, \sigma^2}^{(\beta)}(\xi)
=\frac{\sigma^{-N\{\beta(N-1)+2\}/2}}{C_N^{(\beta)}}
e^{-|\x|^2/2 \sigma^2} |h_N(\x)|^{\beta}, \\
\label{eqn:muba}
&& \mu_{N, \sigma^2}^{(\beta, a)}(\xi)
=\frac{\sigma^{-N \{\beta(N-1)+2(a+1) \} }}{C_N^{(\beta,a)}}
\prod_{j=1}^{N}( x_j^{a} e^{-x_j/2 \sigma^2})
|h_N(\x)|^{\beta}, 
\end{eqnarray}
where 
$$
h_N(\x)=\prod_{1 \leq j < k \leq N} (x_k-x_j)
=\det_{1 \leq j, k \leq N} [x_j^{k-1}],
$$
$|\x|^2=\sum_{j=1}^{N}x_j^2$, and
the normalization factors are given by
\begin{eqnarray}
&& C_N^{(\beta)}=\frac{(2\pi)^{N/2}}{N!} 
\prod_{j=1}^{N} \frac{\Gamma(j \beta/2+1)}{\Gamma(\beta/2+1)},
\nonumber\\
&&
C_N^{(\beta, a)}=\frac{2^{N\{\beta(N-1)+2(a+1)\}/2}}{N!}
\prod_{j=1}^{N} \frac{\Gamma(j \beta/2+1) \Gamma(j\beta/2+a-\beta/2+1)}
{\Gamma(\beta/2+1)}.
\nonumber
\end{eqnarray}
They are the probability density functions
of random configurations
$\Xi=\sum_{j=1}^{N} \delta_{X_j}$ and
$\widetilde{\Xi}=\sum_{j=1}^{N} \delta_{\widetilde{X}_j}$ 
in which the configuration spaces of particle positions 
$\X=(X_1, \dots, X_N)$ and 
$\widetilde{\X}=(\widetilde{X}_1, \dots, \widetilde{X}_N)$ 
are given by
\begin{eqnarray}
&& \X \in \W_N^{\rm A} 
\equiv\{\x=(x_1, \dots, x_N) : x_1 < x_2 < \cdots < x_N\},
\nonumber\\
&&
\widetilde{\X} \in \W_N^{+} 
\equiv\{\x=(x_1, \dots, x_N) : 0 \leq x_1 < x_2 < \cdots < x_N\}.
\nonumber
\end{eqnarray}
In particular, when $\beta=1,2$ and 4,
(\ref{eqn:mub}) gives the distributions of
eigenvalues of $N \times N$ Hermitian random matrices in 
the Gaussian orthogonal ensemble (GOE), GUE,
and the Gaussian symplectic ensemble (GSE)
with variances $\sigma^2$,
respectively \cite{Meh04,For10}.
Similarly, for $\nu \in \N_0$,
(\ref{eqn:muba}) with $(\beta,a)=(1,(\nu-1)/2),
(2, \nu), (4, 2\nu+1)$ give
the distributions of squares of (distinct) singular
values of $(N+\nu) \times N$ random matrices
in the chiral Gaussian orthogonal ensemble (chGOE),
chGUE, and the chiral Gaussian symplectic ensemble (chGSE),
respectively \cite{VZ93,Ver94,SV98}.
Moreover \cite{KT07a},
(i) (\ref{eqn:muba}) with $(\beta,a)=(1, \nu/2)$ was
called `the Laguerre ensemble $\beta=1$ initial condition'
in \cite{FNH99},
(ii) (\ref{eqn:muba}) with $(\beta,a)=(1,0)$ gives the
distribution of squares of positive eigenvalues of
$2N \times 2N$ random matrices in the Gaussian ensemble
of class CI studied by Altland and Zirnbauer \cite{AZ96,AZ97},
(iii) (\ref{eqn:muba}) with $(\beta,a)=(1,-1/2)$ gives the
distribution of squares of positive eigenvalues of
$2N \times 2N$ random matrices in the Gaussian ensemble
of `the real-component version of class D' of the
Bogoliubov-de Gennes universality class \cite{KT04},
and (iv) (\ref{eqn:muba}) with
$(\beta,a)=(4,0)$ and $(4,2)$ give
the distributions of squares of distinct eigenvalues of
$2N \times 2N$ random matrices in the Gaussian ensembles
of class DIII-even \cite{AZ96,AZ97}
and of class DIII-odd \cite{Iva01}, respectively.
We write the expectations of measurable
functions of $\Xi$ and $\widetilde{\Xi}$,
 $F(\Xi)$ and $\widetilde{F}(\widetilde{\Xi})$, 
with respect to distributions 
(\ref{eqn:mub}) and (\ref{eqn:muba}) as
$\bE^{(\beta)}_{N, \sigma^2}[F(\Xi)]$
and $\bE^{(\beta, a)}_{N, \sigma^2}[\widetilde{F}(\widetilde{\Xi})]$,
respectively.
We can say that \cite{FNH99} 
the distributions with $\beta=2$
have {\it unitary symmetry} and 
those with $\beta=1$ do
{\it orthogonal symmetry}.

Recently, we studied the noncolliding BM and 
the noncolliding BESQ$^{(\nu)}, \nu >-1$,
starting not from any fixed configurations
but from the distributions having unitary symmetry, 
$\mu^{(2)}_{N, \sigma^2}$ and $\mu^{(2, \nu)}_{N, \sigma^2}$,
respectively.
We showed that in these cases the determinantal structures
of multitime correlation functions are maintained
but the correlation kernels are replaced by
the time shift $t \to t+\sigma^2$
of the correlation kernels for the special initial
configuration $\xi=N \delta_0$,
{\it i.e.}, the configuration
in which all $N$ particles are put on the origin \cite{K11}.
That is, the equalities
\begin{eqnarray}
&& \rho^{\mu^{(2)}_{N, \sigma^2}}
(t_1, \x^{(1)}_{N_1}; \dots; t_M, \x^{(M)}_{N_M})
\equiv \bE^{(2)}_{N, \sigma^2}
\Big[ \rho^{\Xi}(t_1, \x^{(1)}_{N_1}; \dots;
t_M, \x^{(M)}_{N_M}) \Big]
\nonumber\\
&& \qquad 
=\det_{
\substack{1 \leq j \leq N_{m}, 1 \leq k \leq N_{n} \\ 1 \leq m, n \leq M}
}
[\mbK^{N \delta_0}(t_m+\sigma^2, x_{j}^{(m)}; 
t_n+\sigma^2, x_{k}^{(n)} )],
\nonumber\\
&& \rho_{\nu}^{\mu^{(2,\nu)}_{N, \sigma^2}}
(t_1, \x^{(1)}_{N_1}; \dots; t_M, \x^{(M)}_{N_M})
\equiv \bE^{(2,\nu)}_{N, \sigma^2}
\Big[ \rho^{\Xi}_{\nu}(t_1, \x^{(1)}_{N_1}; \dots;
t_M, \x^{(M)}_{N_M}) \Big]
\nonumber\\
&& \qquad 
=\det_{
\substack{1 \leq j \leq N_{m}, 1 \leq k \leq N_{n} \\ 1 \leq m, n \leq M}
}
[\mbK^{N \delta_0}_{\nu}(t_m+\sigma^2, x_{j}^{(m)}; 
t_n+\sigma^2, x_{k}^{(n)} )], 
\quad \nu > -1,
\label{eqn:shift}
\end{eqnarray}
hold for any $M \in \N$, 
$0 < t_1 < \cdots < t_M < \infty$, $\sigma^2>0$,
$\x^{(m)}_{N_m} \in \W_{N_m}^{\rm A}$ or 
$\x^{(m)}_{N_m} \in \W_{N_m}^{+}$,
$1 \leq N_m \leq N, 1 \leq m \leq M$.
We should note that $\mbK^{N\delta_0}$
and $\mbK^{N\delta_0}_{\nu}$ are 
the correlation kernels known as the 
{\it extended Hermite kernel} 
and {\it extended Laguerre kernel},
respectively \cite{NF98,FNH99,TW04,KT07b,For10}.

In the present paper, we report the cases when
the noncolliding BM and the noncolliding BESQ$^{(\nu)}$
start from the distributions having orthogonal symmetry, 
$\mu^{(1)}_{N, \sigma^2}$ and $\mu^{(1, a)}_{N, \sigma^2}$,
respectively.
For $N \in \N$ and a skew-symmetric $2N \times 2N$
matrix $A=(a_{jk})$, the Pfaffian is defined as
\begin{equation}
\Pf(A) = \Pf_{1 \leq j < k \leq 2N}(a_{jk}) 
= \frac{1}{N !} {\sum_{\pi}}'
{\rm sgn}(\pi) a_{\pi(1) \pi(2)} a_{\pi(3) \pi(4)} 
\cdots a_{\pi(2N-1) \pi(2N)},
\label{eqn:pfaffian}
\end{equation}
where the summation $\sum_{\pi}'$ 
is extended over all permutations $\pi$
of $(1,2,\dots, 2N)$ with restriction
$\pi(2k-1) < \pi(2k), k=1,2,\dots, N$.
The main result of the present paper is the fact that
for any $\sigma^2>0$
we can explicitly determine the $2 \times 2$ {\it skew-symmetric
matrix-valued correlation kernels}
\begin{equation}
\mbA(s,x;t,y;\sigma^2)
=\left( \begin{array}{ll}
\rA_{11}(s,x;t,y;\sigma^2) & \rA_{12}(s,x;t,y;\sigma^2) \cr
-\rA_{12}(t,y;s,x;;\sigma^2) & \rA_{22}(s,x;t,y;;\sigma^2)
\end{array} \right), 
\label{eqn:A1}
\end{equation}
$(x,y) \in \R^2, (s,t) \in [0, \infty)^2$, 
and
\begin{equation}
\mbA^{(\nu,\kappa)}(s,x;t,y;\sigma^2)
=\left( \begin{array}{ll}
\rA_{11}^{(\nu,\kappa)}(s,x;t,y;\sigma^2) 
& \rA_{12}^{(\nu,\kappa)}(s,x;t,y;\sigma^2) \cr
-\rA_{12}^{(\nu,\kappa)}(t,y;s,x;\sigma^2) 
& \rA_{22}^{(\nu,\kappa)}(s,x;t,y;\sigma^2)
\end{array} \right), 
\label{eqn:Anu1}
\end{equation}
$(x.y) \in (0, \infty)^2, (s,t) \in [0, \infty)^2$,
with $\kappa=2(\nu-a)$
such that
\begin{eqnarray}
&& \rho^{\mu^{(1)}_{N, \sigma^2}}
(t_1, \x^{(1)}_{N_1}; \dots; t_M, \x^{(M)}_{N_M})
\equiv \bE^{(1)}_{N, \sigma^2}
\Big[ \rho^{\Xi}(t_1, \x^{(1)}_{N_1}; \dots;
t_M, \x^{(M)}_{N_M}) \Big]
\nonumber\\
\label{eqn:main1}
&& \quad 
=\Pf_{
\substack{1 \leq j \leq N_{m}, 1 \leq k \leq N_{n} \\ 1 \leq m, n \leq M}
}
[\mbA(t_m, x_{j}^{(m)}; 
t_n, x_{k}^{(n)};\sigma^2)],
\\
&& \rho_{\nu}^{\mu^{(1,a)}_{N, \sigma^2}}
(t_1, \x^{(1)}_{N_1}; \dots; t_M, \x^{(M)}_{N_M})
\equiv \bE^{(1,a)}_{N, \sigma^2}
\Big[ \rho^{\Xi}_{\nu}(t_1, \x^{(1)}_{N_1}; \dots;
t_M, \x^{(M)}_{N_M}) \Big]
\nonumber\\
\label{eqn:main2}
&& \quad 
=\Pf_{
\substack{1 \leq j \leq N_{m}, 1 \leq k \leq N_{n} \\ 1 \leq m, n \leq M}
}
[\mbA^{(\nu,\kappa)}(t_m, x_{j}^{(m)}; 
t_n, x_{k}^{(n)};\sigma^2 )], 
\quad \nu > -1, a \in (-1, \nu], 
\end{eqnarray}
hold for any 
$0 < t_1 < \cdots < t_M < \infty$, 
$\x^{(m)}_{N_m} \in \W_{N_m}^{\rm A}$ or 
$\x^{(m)}_{N_m} \in \W_{N_m}^{\rm +}$,
$1 \leq N_m \leq N, 1 \leq m \leq M$.
As an analogue of a determinantal process,
an interacting particle system is said to be
a {\it Pfaffian process}, if any multitime
correlation function is given by a Pfaffian
\cite{MP83,PM83,MM91,MMN98,NF99,FNH99,NKT03,Nag03,KNT04,KT07a,Nag07}.
Then we can state that noncolliding diffusion processes,
which are determinantal processes if they start
from fixed initial configurations,
behave as Pfaffian processes
when they start from distributions
having orthogonal symmetry.

The present paper is organized as follows.
In Section 2 preliminaries and main results
are given.
In Section 3 is devoted to proofs of results.

%%%%%%%%%%%%%%%%%%%%%%%%%%%%%%%%%%%%%%%%%%%%%%%%%%%%%%%%%%
%%%  SEC2   %%%%%%%%%%%%%%%%%%%%%%%%%%%%%%%%%%%%%%%%%%%%%%
%%%%%%%%%%%%%%%%%%%%%%%%%%%%%%%%%%%%%%%%%%%%%%%%%%%%%%%%%%
\SSC{Preliminaries and Main Results}%%%%%%%%%%%%%%%%%%%%%%
%%%%%%%%%%%%%%%%%%%%%%%%%%%%%%%%%%%%%%%%%%%%%%%%%%%%%%%%%%

Let $N \in \N, 0 < T < \infty$, and choose an initial 
configuration $\xi=\sum_{j=1}^{N} \delta_{x_j},
x_1 \leq x_2 \leq \cdots \leq x_N$.
Then consider the $N$-particle system of one-dimensional
standard BMs starting from $\xi$ at time $t=0$
{\it conditioned never to collide with each other}
during time period $(0, T]$ \cite{KT02}.
If the initial configuration is $\xi=N \delta_0$,
that is, all $N$ particles start from
the origin, we can show that the multitime joint
probability density function 
for arbitrary $M+1$ sequence of times
$0 < t_1 < \cdots < t_{M} < t_{M+1} \equiv T$, $M \in \N$, 
is given by the formula
\begin{eqnarray}
&& p_T^{N\delta_0}(t_1, \xi^{(1)}; \dots; t_M, \xi^{(M)}; 
t_{M+1}, \xi^{(M+1)})
= C_{N, T}(t_1) {\rm sgn}(h_N(\x^{(M+1)}))
\nonumber\\
&& \qquad \times
\prod_{m=1}^{M} f(t_{m+1}-t_m, \x^{(m+1)}|\x^{(m)}) 
h_N(\x^{(1)}) 
\prod_{j=1}^{N} p(t_1, x^{(1)}_j|0),
\label{eqn:pN0}
\end{eqnarray}
where $\xi^{(m)}=\sum_{j=1}^{N} \delta_{x^{(m)}_j},
\x^{(m)}=(x^{(m)}_1, \dots, x^{(m)}_N) \in \W_N^{\rm A},
1 \leq m \leq M+1$, 
$$
f(t, \y|\x)=\det_{1 \leq j, k \leq N}
[p(t, y_j|x_k)], \quad
\x, \y \in \W_N^{\rm A}, t \geq 0,
$$
and $C_{N,T}(t)=\pi^{N/2} \{\prod_{j=1}^N \Gamma(j/2)\}^{-1}
T^{N(N-1)/4} t^{-N(N-1)/2}$ \cite{NKT03,KNT04}.
The process, whose finite dimensional distributions
are determined by the formula (\ref{eqn:pN0}),
is temporally inhomogeneous \cite{KT02,KT_Sugaku}. 
In this paper we call it `the noncolliding BM with duration $T$
starting from $N \delta_0$' and
express it by $(\Xi_T(t), t \in [0, T], \P^{N \delta_0})$.

When we take the limit $T \to \infty$, we have a temporally
homogeneous system \cite{KT02,KT07b}, 
which we simply call the noncolliding BM
(starting from $N \delta_0$).
For the noncolliding BM, multitime joint probability density
function is given by the following for
an arbitrary initial configuration $\xi$
with $\xi(\R)=N \in \N$,
\begin{equation}
p^{\xi}(t_1, \xi^{(1)}; \dots; t_M, \xi^{(M)})
=h_N(\x^{(M)}) \prod_{m=1}^{M-1} 
f(t_{m+1}-t_{m}, \x^{(m+1)}|\x^{(m)})
h_N^{(+)}(t_1, \x^{(1)}; \xi),
\label{eqn:pxi}
\end{equation}
$0 < t_1 < \cdots < t_M < \infty$, with
$$
h_N^{(+)}(t, \y; \xi)
=\det_{1 \leq j, k \leq N} \left[
\frac{1}{2 \pi i} \oint_{\rC(\xi_j)} dz \,
\frac{p(t, y_k|z)}{\prod_{x \in \xi_j}(z-x)} \right], 
\quad t \geq 0, \y \in \W_N^{\rm A},
$$
where, for a given initial configuration
$\xi=\sum_{j=1}^{N} \delta_{x_j}, 
x_1 \leq x_2 \leq \cdots \leq x_N$, we define
$\xi_j=\sum_{k=1}^{j} \delta_{x_k}, 1 \leq j \leq N$, 
and $\rC(\xi)$ denotes a closed contour on the
complex plane $\C$ encircling the points in 
$\supp \xi \equiv \{x \in \Lambda: \xi(\{x\}) > 0\}$ 
once in the positive direction \cite{BK05,KT10}.
The noncolliding BM starting from $\xi$ is the temporally
homogeneous process, whose finite dimensional distributions
are determined by (\ref{eqn:pxi}),
and is denoted by $(\Xi(t), t \in [0, \infty), \P^{\xi})$
in this paper.
We can prove that 
$\Xi(t, \cdot)=\sum_{j=1}^{N} \delta_{X_j(t)}(\cdot), t \geq 0$
solves the following system of stochastic differential
equations (SDEs),
$$
dX_j(t)=dB_j(t)+ 
\sum_{1 \leq k \leq N, k \not= j}
\frac{dt}{X_j(t)-X_k(t)},
\quad
1 \leq j \leq N, \quad t \geq 0,
$$
with independent one-dimensional standard
BMs $\{B_j(t) \}_{j=1}^{N}$,
which is the $\beta=2$ case of Dyson's BM model
\cite{KT07b,KT_Sugaku}.

In \cite{KT07a}, a temporally inhomogeneous
noncolliding diffusion process was introduced,
which is called {\it the noncolliding squared generalized meander}
with duration $T$.
It is a two-parameter family of processes
indexed by $\nu > -1$ and
$\kappa \in [0, 2(\nu+1))$ starting from
the configuration $N \delta_0, N \in \N$,
which includes the processes studied in
\cite{FNH99} and \cite{Nag03} 
as special cases.
This family of processes is denoted here by
$(\Xi_T^{(\nu,\kappa)}(t), t \in [0, T], \P^{N\delta_0})$.
The multitime joint probability density function is given by
\begin{eqnarray}
&& p_{T,(\nu,\kappa)}^{N\delta_0}
 (t_1, \xi^{(1)}; \dots; t_M, \xi^{(M)}; t_{M+1}, \xi^{(M+1)})
=C_{N, T}^{(\nu,\kappa)}(t_1) {\rm sgn}(h_N(\x^{(M+1)}))
\nonumber\\
&& \qquad \times
\prod_{m=1}^{M} f^{(\nu,\kappa)}(t_{m+1}-t_m, \x^{(m+1)}|\x^{(m)})
h_N(\x^{(1)}) \prod_{j=1}^{N} p^{(\nu,\kappa)}(t_1, x^{(1)}_j|0),
\label{eqn:pN0nk}
\end{eqnarray}
$0 < t_1 < \cdots < t_{M} < t_{M+1} \equiv T$, where
\begin{eqnarray}
p^{(\nu,\kappa)}(t, y|x)
&=&  \left\{ \begin{array}{ll}
\displaystyle{
\frac{1}{2t} \left( \frac{y}{x} \right)^{(\nu-\kappa)/2}
e^{-(x+y)/2t}
I_{\nu} \left( \frac{\sqrt{xy}}{t} \right)},
& \quad t >0, x >0, y \geq 0, \cr
\displaystyle{
\frac{y^{\nu-\kappa/2}e^{-y/2t} }{(2 t)^{\nu+1} \Gamma(\nu+1)} },
& \quad t >0, x=0, y \geq 0, \cr
& \cr
\delta(y-x),
& \quad t=0, x,y \geq 0,
\end{array} \right.
\label{eqn:pnk}
\end{eqnarray}
$$
f^{(\nu,\kappa)}(t, \y|\x)=\det_{1 \leq j, k \leq N}
[p^{(\nu, \kappa)}(t, y_j|x_k)], \quad
\x, \y \in \W_N^{+}, t \geq 0,
$$
and 
$$
C_{N,T}^{(\nu,\kappa)}(t)
=\frac{T^{(N+\kappa-1)N/2} t^{-(N-1)N}}{2^{N(N-\kappa-1)/2}}
\prod_{j=1}^{N} \frac{\Gamma(\nu+1) \Gamma(1/2)}
{\Gamma \left(j/2 \right) \Gamma \left( (j+1+2\nu-\kappa)/2 \right)}.
$$
On the other hand, for the noncolliding BESQ$^{(\nu)}$,
which is obtained as temporally homogeneous limit $T \to \infty$
of $\Xi_T^{(\nu,\kappa)}(t)$ \cite{KT04,KT07a},
the multitime joint probability density function
is obtained for an arbitrary initial configuration
$\xi=\sum_{j=1}^{N} \delta_{x_j},
0 \leq x_1 \leq x_2 \leq \cdots \leq x_N, N \in \N$ as
\begin{equation}
p^{\xi}_{\nu}(t_1, \xi^{(1)}; \dots; t_M, \xi^{(M)})
=h_N(\x^{(M)}) \prod_{m=1}^{M-1} 
f^{(\nu)}(t_{m+1}-t_{m}, \x^{(m+1)}|\x^{(m)})
h_N^{(\nu,+)}(t_1, \x^{(1)}; \xi),
\label{eqn:pxin}
\end{equation}
$0 < t_1 < \cdots < t_M < \infty$, with \cite{DF08,KT11}
\begin{eqnarray}
&& f^{(\nu)}(t, \y|\x)
=\det_{1 \leq j, k \leq N} 
[ p^{(\nu)}(t, y_j|x_k)],
\quad \x, \y \in \W_N^{+}, t \geq 0,
\nonumber\\
&&
h_N^{(\nu,+)}(t, \y; \xi)
=\det_{1 \leq j, k \leq N} \left[
\frac{1}{2 \pi i} \oint_{\rC(\xi_j)} dz \,
\frac{p^{(\nu)}(t, y_k|z)}{\prod_{x \in \xi_j}(z-x)} \right].
\nonumber
\end{eqnarray}
We write the noncolliding BESQ$^{(\nu)}$ starting from $\xi$
as $(\Xi^{(\nu)}(t), t \in [0, \infty), \P^{\xi})$.
If we set $\Xi^{(\nu)}(t, \cdot)=\sum_{j=1}^{N}
\delta_{\widetilde{X}_j(t)}(\cdot)$,
$\widetilde{\X}(t)=(\widetilde{X}_1(t), \dots,
\widetilde{X}_N(t))$ satisfies the SDEs
\begin{eqnarray}
d \widetilde{X}_j(t) &=& 2 \sqrt{\widetilde{X}_j(t)} 
d \widetilde{B}_j(t) + 2 (\nu+1) dt
\nonumber\\
&& + 4 \widetilde{X}_j(t)
\sum_{1 \leq k \leq N, k \not= j}
\frac{dt}{\widetilde{X}_j(t)-\widetilde{X}_k(t)},
\quad 1 \leq j \leq N, \quad t \geq 0,
\nonumber
\end{eqnarray}
where $\{\widetilde{B}_j(t)\}_{j=1}^{N}$
are independent one-dimensional standard BMs
and, if $-1 < \nu < 0$,
the reflection boundary condition is assumed at
the origin \cite{KT11}.

We write the noncolliding BM starting from the GOE eigenvalue
distribution $\mu^{(1)}_{N, \sigma^2}$ as
$(\Xi(t), t \in [0, \infty), \P^{\mu^{(1)}_{N, \sigma^2}})$
and the noncolliding BESQ$^{(\nu)}, \nu>-1$ from the
distribution having orthogonal symmetry,
$\mu^{(1,a)}_{N, \sigma^2}, 
a \in (-1, \nu]$,
as $(\Xi^{(\nu)}(t), t \in [0, \infty), 
\P^{\mu^{(1,a)}_{N, \sigma^2}})$.

In general, two processes having the same state space
are said to be {\it equivalent} if they have the same
finite-dimensional distributions, that is,
if, for any finite sequence of times 
$0 < t_1 < \cdots < t_M < \infty, M \in \N$,
the multitime joint probability density functions
coincide with each other \cite{RY98}.
The key lemma of the present study is the
following equivalence.

%%%%%%%%%%%%%%%%%%%%%%%%%%%%%%%%%
\begin{lem}
\label{thm:equivalence}
For $\sigma^2>0$, let
\begin{equation}
c_{\sigma^2}(t)=\frac{\sigma^2}{\sigma^2+t},
\quad t \in [0, \infty).
\label{eqn:cs}
\end{equation}
Then
\begin{equation}
(\Xi(t), t \in [0, \infty), \P^{\mu^{(1)}_{N, \sigma^2}})
= \left( \frac{1}{c_{\sigma^2}(t)} \circ
\Xi_{\sigma^2}(\sigma^2 c_{\sigma^2}(t)), t \in [0, \infty),
\P^{N \delta_0} \right).
\label{eqn:Eq1}
\end{equation}
If the relation
\begin{equation}
a=\nu-\frac{\kappa}{2},
\quad \nu > -1, \quad
a \in(-1, \nu], 
\label{eqn:ank}
\end{equation}
is satisfied, then
\begin{equation}
(\Xi^{(\nu)}(t), t \in [0, \infty), \P^{\mu^{(1,a)}_{N, \sigma^2}})
= \left( \frac{1}{c_{\sigma^2}(t)^2} \circ
\Xi^{(\nu, \kappa)}_{\sigma^2}(\sigma^2 c_{\sigma^2}(t)), 
t \in [0, \infty),
\P^{N \delta_0} \right).
\label{eqn:Eq2}
\end{equation}
\end{lem}
\vskip 0.5cm
%%%%%%%%%%%%%%%%%%%%%%%%%%%%%%%%%%

Remark that
\begin{equation}
0 < s < t < \infty
\quad \Longleftrightarrow \quad
\sigma^2 > \sigma^2 c_{\sigma^2}(s)
> \sigma^2 c_{\sigma^2}(t) > 0.
\label{eqn:reverse}
\end{equation}
Therefore, the RHS of (\ref{eqn:Eq1}) and (\ref{eqn:Eq2})
are time reverses of the processes
$\Xi_{T}(\cdot)$ and $\Xi^{(\nu,\kappa)}_{T}(\cdot)$
with duration $T=\sigma^2$,
followed by dilatation with factors 
$1/c_{\sigma^2}(\cdot)$ and $1/c_{\sigma^2}(\cdot)^2$,
respectively.

In \cite{NKT03,KNT04} and \cite{KT07a}, it was proved that
$(\Xi_T(t), t \in [0, T], \P^{N \delta_0})$ 
and $(\Xi_T^{(\nu, \kappa)}(t), t \in [0, T], \P^{N \delta_0})$,
$0 < T < \infty$, $\nu>-1$, $\kappa \in [0, 2(\nu+1))$
are Pfaffian processes, respectively. 
Then by the equivalence (\ref{eqn:Eq1}) and (\ref{eqn:Eq2})
of Lemma \ref{thm:equivalence},
the following main Theorems are obtained.
(For simplicity of expressions, we show the elements
of the matrix-valued correlation kernels 
only for the case that the number of particles $N$ is even.
See, for example, \cite{Nag07}
for the general theory of Pfaffian expressions
of correlation functions.)
Let $H_n$ and $L_n^{\nu}$ be the Hermite polynomial of degree $n$
and the Laguerre polynomials of degree $n$ with index $\nu$;
\begin{eqnarray}
H_n(x) &=& n ! \sum_{k=0}^{[n/2]}
(-1)^k \frac{(2x)^{n-2k}}{k ! (n-2k)!},
\nonumber\\
L_n^{\nu}(x) &=& \sum_{k=0}^{n} (-1)^k
\frac{\Gamma(n+\nu+1) x^k}
{\Gamma(k+\nu+1) (n-k)! k!},
\quad n \in \N_0,
\nonumber
\end{eqnarray}
where $[r]$ denotes the greatest integer not 
greater than $r$.
For $n \in \Z \equiv \{\dots, -1,0,1,2, \dots\}$
and $\alpha \in \R$ we define
$$
{n+\alpha \choose n} =
\left\{
   \begin{array}{ll}
 \displaystyle{\frac{\Gamma (n+\alpha+1)}{\Gamma(n+1)\Gamma(\alpha+1)},}
 & \mbox{if $n\in\N$, $\alpha \notin \Z_-$}, \\
 \displaystyle{\frac{(-1)^n \Gamma (-\alpha)}{\Gamma(n+1)\Gamma(-n-\alpha)},}
 & \mbox{if $n\in\N$, $n+\alpha\in \Z_-$}, \\
    0,
 & \mbox{if $n\in\N$, $\alpha \in \Z_-$, $n+\alpha\in \N_0$}, \\
   1,
 & \mbox{if $n=0$}, \\
   0,
 & \mbox{if $n\in\Z_-$}, \\
   \end{array}\right.
$$
where $\Z_{-}=\Z \setminus \N_0$.

%%%%%%%%%%%%%%%%%%%%%%%%%%%%%%%%%
\begin{thm}
\label{thm:Main1}
The noncolliding BM with a finite number of particles
$N \in \N$, starting from the GOE eigenvalue distribution
with variance $\sigma^2>0$,
$(\Xi(t), t \in [0, \infty), \P^{\mu^{(1)}_{N, \sigma^2}})$,
is a Pfaffian process.
When $N$ is even, the elements of the matrix-valued
correlation kernel (\ref{eqn:A1}) are given by
\begin{eqnarray}
&& A_{11}(s,x;t,y;\sigma^2)
=\sum_{k=0}^{N/2-1} \frac{1}{d_k(\sigma^2)}
\Big[ B_{2k}(s,x;\sigma^2) B_{2k+1}(t,y; \sigma^2)
-B_{2k+1}(s,x; \sigma^2) B_{2k}(t,y; \sigma^2) \Big],
\nonumber\\
&& A_{12}(s,x;t,y;\sigma^2) \nonumber\\
&& \quad = \left\{ \begin{array}{ll}
\displaystyle{\sum_{k=0}^{N/2-1} \frac{1}{d_{k}(\sigma^2)}
\Big[ B_{2k+1}(s,x; \sigma^2) C_{2k}(t,y;\sigma^2) 
-B_{2k}(s,x; \sigma^2) C_{2k+1}(t,y;\sigma^2) \Big]}, 
& \mbox{if $s \leq t$}, \cr
\displaystyle{-\sum_{k=N/2}^{\infty} \frac{1}{d_{k}(\sigma^2)}
\Big[ B_{2k+1}(s,x; \sigma^2) C_{2k}(t,y;\sigma^2)
-B_{2k}(s,x; \sigma^2) C_{2k+1}(t,y;\sigma^2) \Big]}, 
& \mbox{if $s > t$},
\end{array} \right.
\nonumber\\
&& A_{22}(s,x;t,y;\sigma^2)
= \sum_{k=N/2}^{\infty} \frac{1}{d_{k}(\sigma^2)}
\Big[ C_{2k}(s,x; \sigma^2) C_{2k+1}(t,y;\sigma^2)
-C_{2k+1}(s,x; \sigma^2) C_{2k}(t,y; \sigma^2) \Big], 
\nonumber\\
\label{eqn:As}
\end{eqnarray}
with
\begin{eqnarray}
d_k(\sigma^2) &=& 2 \sigma^2 \Gamma(k+1/2) \Gamma(k+1),
\nonumber\\
B_{2k}(s,x; \sigma^2)
&=& \left(\frac{\sigma^2+2s}{4 \sigma^2}\right)^k
e^{-x^2/2(\sigma^2+s)} H_{2k}
\left( \frac{x}{\sqrt{\sigma^2+2s}} \right),
\nonumber\\
B_{2k+1}(s,x;\sigma^2) 
&=& \left(\frac{\sigma^2+2s}{4 \sigma^2}\right)^{(2k+1)/2}
e^{-x^2/2(\sigma^2+s)} 
\nonumber\\
&\times& 
\left\{ H_{2k+1} \left( \frac{x}{\sqrt{\sigma^2+2s}} \right)
- \frac{4 k\sigma^2}{\sigma^2+2s}
H_{2k-1} \left( \frac{x}{\sqrt{\sigma^2+2s}} \right) \right\},
\nonumber\\
C_{2k}(s,x; \sigma^2)
&=&
\frac{(2k)!}{k!} 2^{-2k+1} \sigma
\sqrt{\frac{\sigma^2}{\sigma^2+2s}}
e^{x^2/2(\sigma^2+s)-x^2/(\sigma^2+2s)}
\nonumber\\
&\times& 
\sum_{\ell=k}^{\infty} \frac{\ell!}{(2\ell+1)!}
\left( \frac{\sigma^2}{\sigma^2+2s} \right)^{(2\ell+1)/2}
H_{2 \ell+1} \left( \frac{x}{\sqrt{ \sigma^2+2s}} \right),
\nonumber\\
C_{2k+1}(s,x; \sigma^2)
&=& - 2^{-2k+1} \sigma \sqrt{\frac{\sigma^2}{\sigma^2+2s}}
\left( \frac{\sigma^2}{\sigma^2+2s} \right)^{k}
\nonumber\\
&\times& 
e^{x^2/2(\sigma^2+s)-x^2/(\sigma^2+2s)}
H_{2k}\left(\frac{x}{\sqrt{\sigma^2+2s}} \right),
\quad k \in \N_0. 
\label{eqn:BC}
\end{eqnarray}
\end{thm}
\vskip 0.5cm
%%%%%%%%%%%%%%%%%%%%%%%%%%%%%%%%%%

%%%%%%%%%%%%%%%%%%%%%%%%%%%%%%%%%
\begin{thm}
\label{thm:Main2}
The noncolliding BESQ$^{(\nu)}, \nu>-1$ 
with a finite number of particles
$N \in \N$, starting from the distribution
having orthogonal symmetry, 
$\mu^{(1,a)}_{N, \sigma^2}, \sigma^2>0, a \in (-1, \nu]$,
$(\Xi^{(\nu)}(t), t \in [0, \infty), 
\P^{\mu^{(1,a)}_{N, \sigma^2}})$,
is a Pfaffian process.
Let $\kappa=2(\nu-a)$. Then,
when $N$ is even, the elements of the matrix-valued
correlation kernel (\ref{eqn:Anu1}) are given by
\begin{eqnarray}
&& A^{(\nu,\kappa)}_{11}(s,x;t,y;\sigma^2)
\nonumber\\
&& \qquad
=\sum_{k=0}^{N/2-1} \frac{1}{d^{(\nu,\kappa)}_k(\sigma^2)}
\Big[ B^{(\nu,\kappa)}_{2k}(s,x;\sigma^2) 
B^{(\nu,\kappa)}_{2k+1}(t,y; \sigma^2)
-B^{(\nu,\kappa)}_{2k+1}(s,x; \sigma^2) 
B^{(\nu,\kappa)}_{2k}(t,y; \sigma^2) \Big],
\nonumber\\
&& A^{(\nu,\kappa)}_{12}(s,x;t,y;\sigma^2) \nonumber\\
&& \qquad =\sum_{k=0}^{N/2-1} \frac{1}{d_{k}(\sigma^2)}
\Big[ 
B^{(\nu,\kappa)}_{2k+1}(s,x; \sigma^2)
C^{(\nu,\kappa)}_{2k}(t,y;\sigma^2) 
-B^{(\nu,\kappa)}_{2k}(s,x; \sigma^2) 
C^{(\nu,\kappa)}_{2k+1}(t,y; \sigma^2) \Big]
\nonumber\\
&& \qquad \quad
-\1(s > t) p^{(\nu,\kappa)}_{-}(s,x;t,y;\sigma^2),
\nonumber\\
&& A^{(\nu,\kappa)}_{22}(s,x;t,y;\sigma^2)
\nonumber\\
&& \qquad 
= \sum_{k=N/2}^{\infty} \frac{1}{d^{(\nu,\kappa)}_{k}(\sigma^2)}
\Big[ C^{(\nu,\kappa)}_{2k}(s,x; \sigma^2) 
C^{(\nu,\kappa)}_{2k+1}(t,y;\sigma^2)
-C^{(\nu,\kappa)}_{2k+1}(s,x; \sigma^2) 
C^{(\nu,\kappa)}_{2k}(t,y; \sigma^2) \Big], 
\nonumber\\
\label{eqn:Anus}
\end{eqnarray}
with
\begin{eqnarray}
d_k^{(\nu,\kappa)}(\sigma^2) &=& 
2^{-2\nu} \sigma^{-2\kappa}
\frac{(2k)! \Gamma(2k+2+2 \nu-\kappa)}
{\Gamma(\nu+1)^2},
\nonumber\\
B^{(\nu,\kappa)}_{k}(s,x; \sigma^2)
&=&  
\frac{k!}{2^{\nu+1} \Gamma(\nu+1)}
\frac{1}{\sigma^{2(\kappa+1)}(\sigma^2+s)^{\nu-\kappa}}
\nonumber\\
&\times&
e^{-x/2(\sigma^2+s)} x^{\nu-\kappa/2}
\sum_{j=0}^{k} \alpha_{k,j}
\left(\frac{\sigma^2+2s}{\sigma^2} \right)^j
L^{\nu}_j \left(\frac{x}{\sigma^2+2s} \right),
\nonumber\\
C^{(\nu,\kappa)}_{2k}(s,x; \sigma^2)
&=&
\frac{(2k)! \Gamma(2\nu-\kappa+1)}{2^{\nu-1} \Gamma(\nu+1)}
\sigma^2 \frac{(\sigma^2+s)^{\nu-\kappa}}{(\sigma^2+2s)^{\nu+1}}
{2k+2\nu-\kappa+1 \choose 2k+1} 
\nonumber\\ 
&\times& e^{x/2(\sigma^2+s)-x/(\sigma^2+2s)} x^{\kappa/2}
\sum_{j=2k+1}^{\infty} \beta_{j, 2k+1}
\frac{\Gamma(j+1)}{\Gamma(j+1+\nu)}
\left( \frac{\sigma^2}{\sigma^2+2s} \right)^j
L^{\nu}_{j} \left(\frac{x}{\sigma^2+2s} \right), 
\nonumber\\
C^{(\nu,\kappa)}_{2k+1}(s,x; \sigma^2)
&=& -
\frac{(2k+1)! \Gamma(2\nu-\kappa+1)}{2^{\nu-1} \Gamma(\nu+1)}
\sigma^2 \frac{(\sigma^2+s)^{\nu-\kappa}}{(\sigma^2+2s)^{\nu+1}}
{2k+2\nu-\kappa+1 \choose 2k+1} 
\nonumber\\ 
&\times& e^{x/2(\sigma^2+s)-x/(\sigma^2+2s)} x^{\kappa/2}
\sum_{j=2k}^{\infty} \beta_{j, 2k}
\frac{\Gamma(j+1)}{\Gamma(j+1+\nu)}
\left( \frac{\sigma^2}{\sigma^2+2s} \right)^j
L^{\nu}_{j} \left(\frac{x}{\sigma^2+2s} \right), 
\nonumber\\
\label{eqn:BCnu}
\end{eqnarray}
$k \in \N_0$, where
\begin{eqnarray}
&& \alpha_{k,j}=\left\{ \begin{array}{ll}
\displaystyle{{k-j+\nu-\kappa \choose k-j}},
& \mbox{if $k$ is even}, \cr
\displaystyle{\frac{k+2\nu-\kappa}{k}
{k-2-j+\nu-\kappa \choose k-2-j}
-{k-j+\nu-\kappa \choose k-j}},
& \mbox{if $k$ is odd},
\end{array} \right.
\nonumber\\
&& \beta_{j,2k}=
{j-2k-\nu+\kappa-2 \choose j-2k},
\quad j \geq 2k,
\nonumber\\
&& \beta_{j,2k+1}
=-\sum_{\ell=k+1}^{[(j+1)/2]}
b(2k+3,2 \ell-1)
{j-2 \ell-\nu+\kappa-1 \choose j-2 \ell+1},
\quad j \geq 2k+1,
\nonumber
\end{eqnarray}
with
$$
b(m,n)=\left\{ \begin{array}{ll}
\displaystyle{\prod_{\ell=0}^{(n-m)/2}
\frac{m+2\ell+2\nu-\kappa}{m+2\ell}},
& \mbox{if $m,m$ are odd and $m \leq n$}, \cr
1, & \mbox{if $m,n$ are odd and $m>n$}, \cr
0, & \mbox{otherwise},
\end{array} \right.
$$
and
$$
p^{(\nu,\kappa)}_{-}(s, x; t, y; \sigma^2)
= \left\{ \begin{array}{ll}
\displaystyle{
\left( \frac{\sigma^2+s}{\sigma^2+t} \right)^{\nu-\kappa}
e^{x/2(\sigma^2+s)-y/2(\sigma^2+t)}
p^{(\nu,\kappa)}(s-t, y|x)},
&x > 0, \cr
\displaystyle{ 
\left(\frac{\sigma^2+t}{\sigma^2} \right)^{-(\nu-\kappa)}
\left( \frac{\sigma^2+s}{\sigma^2} \right)^{\kappa}
e^{-y/2(\sigma^2+t)} p^{(\nu,\kappa)}(s-t,y|0)},
&x=0,
\end{array} \right.
$$
for $s>t, y \geq 0$.
\end{thm}
\vskip 0.5cm
%%%%%%%%%%%%%%%%%%%%%%%%%%%%%%%%%%

%%%%%%%%%%%%%%%%%%%%%%%%%%%%%%%%%%%%%%%%%%%%%%%%%%%%%%%%%%
%%%  SEC3 %%%%%%%%%%%%%%%%%%%%%%%%%%%%%%%%%%%%%%%%%%%
%%%%%%%%%%%%%%%%%%%%%%%%%%%%%%%%%%%%%%%%%%%%%%%%%%%%%%%%%%
\SSC{Proofs of Theorems}%%%%%%%%
%%%%%%%%%%%%%%%%%%%%%%%%%%%%%%%%%%%%%%%%%%%%%%%%%%%%%%%%%%
\subsection{Proof of Lemma \ref{thm:equivalence}}
%%%%%%%%%%%%%%%%%%%%%%%%%%%%%%%%%%%%%%%%%%%%%%%%%%%%%

When the initial configuration $\xi=\sum_{j=1}^{N} \delta_{x_j}$
is distributed according to $\mu^{(1)}_{N, \sigma^2}, \sigma^2>0$,
there is no multiple point in $\x=(x_1, \dots, x_N)$
with probability one, {\it i.e.},
$\bP^{(1)}_{N, \sigma^2}[\X \in \W_N^{\rm A}]=1$.
In this case
$
h_N^{(+)}(t, \y; \xi)=f(t, \y|\x)/h_N(\x),
$
and we can confirm the equality
\begin{eqnarray}
&& h_N(\x^{(M)}) h^{(+)}_N(t_1, \x^{(1)}; \xi)
\mu^{(1)}_{N, \sigma^2}(\xi)
\nonumber\\
&& \quad = C_{N, \sigma^2}(\sigma^2 c_{\sigma^2}(t_M))
h_N(c_{\sigma^2}(t_M) \x^{(M)})
\prod_{j=1}^{N} p(\sigma^2 c_{\sigma^2}(t_M), 
c_{\sigma^2}(t_M) x^{(M)}_j | 0)
\nonumber\\
&& \qquad \qquad \times
c_{\sigma^2}(t_M)^{N/2} c_{\sigma^2}(t_1)^{N/2}
e^{|\x^{(M)}|^2/2(\sigma^2+t_M)-|\x^{(1)}|^2/2(\sigma^2+t_1)}
\nonumber\\
&& \qquad \qquad \times
f(\sigma^2-\sigma^2 c_{\sigma^2}(t_1), \x|
c_{\sigma^2}(t_1) \x^{(1)})
{\rm sgn}(h_N(\x)),
\nonumber
\end{eqnarray}
where, for $c>0, \y=(y_1, \dots, y_N) \in \W_N^{\rm A}$,
we put $c \y \equiv (cy_1, \dots, c y_N)$.
Similarly, we can see the equalities
for $1 \leq m \leq M-1$,
\begin{eqnarray}
&& f(t_{m+1}-t_{m}, \x^{(m+1)}|\x^{(m)})
\nonumber\\
&& \qquad = c_{\sigma^2}(t_m)^{N/2} 
c_{\sigma^2}(t_{m+1})^{N/2} 
e^{|\x^{(m)}|^2/2(\sigma^2+t_m)-|\x^{(m+1)}|^2/2(\sigma^2+t_{m+1})}
\nonumber\\
&& \qquad \times
f(\sigma^2 c_{\sigma^2}(t_m)-\sigma^2 c_{\sigma^2}(t_{m+1}),
c_{\sigma^2}(t_m) \x^{(m)}| c_{\sigma^2}(t_{m+1}) \x^{(m+1)}).
\nonumber
\end{eqnarray}
Then for any $M \in \N$, 
$0 < t_1 < t_2 < \cdots < t_M < \infty$,
$\xi=\sum_{j=1}^{N} \delta_{x_j},
\x \in \W_N^{\rm A}$, and 
$\xi^{(m)}=\sum_{j=1}^{N} \delta_{x^{(m)}_j},
\x^{(m)} \in \W_N^{\rm A}, 1 \leq m \leq M$,
the equality 
\begin{eqnarray}
&& \mu^{(1)}_{N, \sigma^2}(\xi)
p^{\xi}(t_1, \xi^{(1)}; \dots ; t_M, \xi^{(M)})
\nonumber\\
&&
= p^{N \delta_0}_{\sigma^2}
(\sigma^2 c_{\sigma^2}(t_M), c_{\sigma^2}(t_M) \circ \xi^{(M)};
\dots; \sigma^2 c_{\sigma^2}(t_1), c_{\sigma^2}(t_1) \circ \xi^{(1)};
\sigma^2, \x)
\prod_{m=1}^{M} c_{\sigma^2}(t_m)^N 
\nonumber\\
\label{eqn:idenA3}
\end{eqnarray}
holds.
Integration of the LHS of (\ref{eqn:idenA3}) over
$\x \in \W_N^{\rm A}$ gives
\begin{eqnarray}
\int_{\W_N^{\rm A}} d \x \,
\mu^{(1)}_{N, \sigma^2}(\xi)
p^{\xi}(t_1, \xi^{(1)}; \dots; t_M, \xi^{(M)})
&=& \bE^{(1)}_{N, \sigma^2}
\Big[ p^{\Xi}(t_1, \xi^{(1)}; \dots; t_M, \xi^{(M)}) \Big]
\nonumber\\
&=& p^{\mu^{(1)}_{N, \sigma^2}}
(t_1, \xi^{(1)}; \dots; t_M, \xi^{(M)}),
\nonumber
\end{eqnarray}
and that of the RHS of (\ref{eqn:idenA3}) gives
the multitime joint probability density function
of the noncolliding BM
starting from $N \delta_0$ with duration $T=\sigma^2$,
in which observations are performed at $M$ times
in the reversed order 
$0< \sigma^2 c_{\sigma^2}(t_M) <\sigma^2 c_{\sigma^2}(t_{M-1}) 
< \cdots < \sigma^2 c_{\sigma^2}(t_1) < \sigma^2$,
multiplied by the scale factors
$c_{\sigma^2}(t_m), 1 \leq m \leq M$.
Then the equivalence of the processes
(\ref{eqn:Eq1}) is concluded. 
In a similar way, we can prove (\ref{eqn:Eq2}).
\qed

%%%%%%%%%%%%%%%%%%%%%%%%%%%%%%%%%%%%%%%%%%%%%%%%%%%%%
\subsection{Proof of Theorem \ref{thm:Main1}}
%%%%%%%%%%%%%%%%%%%%%%%%%%%%%%%%%%%%%%%%%%%%%%%%%%%%%

For a sequence $(N_m)_{m=1}^{M}$ of positive integers
less than or equal to $N$,
the $(N_1, \dots, N_M)$-multitime correlation function
at $M$ times $0<t_1< \cdots < t_M < T$ of
$(\Xi_T(t), t \in [0, T], \P^{N \delta_0})$
is obtained from (\ref{eqn:pN0}) by
\begin{eqnarray}
&& \rho^{N \delta_0}_T (t_{1}, \x^{(1)}_{N_1} ; 
\dots; t_M, \x^{(M)}_{N_M}) 
\nonumber\\
&&=\prod_{m=1}^{M} \int_{\R^{N-N_{m}}}
\prod_{j=N_m+1}^{N}
\frac{d x^{(m)}_j}{(N-N_m)!}
\int_{\R^N} \frac{d \x^{(M+1)}}{N!}
p^{N \delta_0}_T (t_1, \xi^{(1)}; 
\dots; t_M, \xi^{(M)}; t_{M+1}, \xi^{(M+1)}).
\nonumber
\end{eqnarray}
In \cite{NKT03,KNT04}, the functions
$\widetilde{\rA}_{11}(s,x;t,y;T,t_1), \widetilde{\rA}_{12}(s,x;t,y;T,t_1),
\widetilde{\rA}_{22}(s,x;t,y;T,t_1), 0 <s,t<T, (x,y) \in \R^2$
are given such that
\begin{eqnarray}
&&\rho^{N \delta_0}_T (t_{1}, \x^{(1)}_{N_1} ; 
\dots; t_M, \x^{(M)}_{N_M}) 
\nonumber\\
&& \qquad
=\Pf_{
\substack{1 \leq j \leq N_{m}, 1 \leq k \leq N_{n} \\ 1 \leq m, n \leq M}
}
\left[ \left(
\begin{array}{ll}
\widetilde{\rA}_{11}(t_m,x^{(m)}_j; t_n, x^{(n)}_k;T,t_1) 
& \widetilde{\rA}_{12}(t_m, x^{(m)}_j;t_n, x^{(n)}_k;T,t_1) \cr
-\widetilde{\rA}_{12}(t_n, x^{(n)}_k;t_m, x^{(m)}_j;T,t_1) 
& \widetilde{\rA}_{22}(t_m, x^{(m)}_j;t_n, x^{(n)}_k;T,t_1)
\end{array} \right) \right].
\nonumber
\end{eqnarray}
By direct calculation we have found that
\begin{eqnarray}
&& \widetilde{\rA}_{jk}(\sigma^2 c_{\sigma^2}(s), c_{\sigma^2}(s)x;
\sigma^2 c_{\sigma^2}(t),c_{\sigma^2}(t)y; \sigma^2, c_{\sigma^2}(t_M))
\nonumber\\
&& \quad
= c_{\sigma^2}(s)^{-1/2} c_{\sigma^2}(t)^{-1/2}
\rA_{jk}(s,x;t,y;\sigma^2),
\label{eqn:AA}
\end{eqnarray}
for $(j,k)=(1,1), (1,2), (2,2)$,
where $\rA_{jk}, (j,k)=(1,1), (1,2), (2,2)$, 
are given by (\ref{eqn:As}) with (\ref{eqn:BC}).
Note that by definition of Pfaffian (\ref{eqn:pfaffian}),
with any set of factors $v_j, 1 \leq j \leq 2N$,
\begin{equation}
\Pf_{1 \leq j < k \leq 2N}(v_j a_{jk} v_k)
=\prod_{j=1}^{2N} v_j \times \Pf_{1 \leq j < k \leq 2N} (a_{jk}).
\label{eqn:Pfaffian2}
\end{equation}
Then by the equality (\ref{eqn:Eq1}) of Lemma \ref{thm:equivalence}
(see also (\ref{eqn:idenA3})),
\begin{eqnarray}
&& \rho^{\mu^{(1)}_{N,\sigma^2}}
(t_1, \x^{(1)}_{N_1}; \dots ; t_M, \x^{(M)}_{N_M})
\nonumber\\
&&
= \Pf_{
\substack{1 \leq j \leq N_{m}, 1 \leq k \leq N_{n} \\ 1 \leq m, n \leq M}
} 
\left[ c_{\sigma^2}(t_m)^{-1}
\begin{footnotesize}
\left( \begin{array}{ll}
\rA_{11}(t_m,x^{(m)}_j; t_n, x^{(n)}_k; \sigma^2) &
\rA_{12}(t_m,x^{(m)}_j; t_n, x^{(n)}_k; \sigma^2) \cr
-\rA_{12}(t_n, x^{(n)}_k; t_m,x^{(m)}_j; \sigma^2) &
\rA_{22}(t_m,x^{(m)}_j; t_n, x^{(n)}_k; \sigma^2)
\end{array} \right)
\end{footnotesize}
c_{\sigma^2}(t_n)^{-1} \right]
\nonumber\\
&& \qquad \qquad \qquad \qquad \times
\prod_{\ell=1}^{M} c_{\sigma^2}(t_{\ell})^N
\nonumber\\
&&= \Pf_{
\substack{1 \leq j \leq N_{m}, 1 \leq k \leq N_{n} \\ 1 \leq m, n \leq M}
} [\mbA(t_m, x_{j}^{(m)}; 
t_n, x_{k}^{(n)};\sigma^2)],
\nonumber
\end{eqnarray}
and the proof is completed. \qed

%%%%%%%%%%%%%%%%%%%%%%%%%%%%%%%%%%%%%%%%%%%%%%%%%%%%%
\subsection{Proof of Theorem \ref{thm:Main2}}
%%%%%%%%%%%%%%%%%%%%%%%%%%%%%%%%%%%%%%%%%%%%%%%%%%%%%

In \cite{KT07a}, the functions
$\widetilde{\rA}^{(\nu,\kappa)}_{11}(s,x;t,y;T,t_1)$, 
$\widetilde{\rA}^{(\nu,\kappa)}_{12}(s,x;t,y;T,t_1)$,
$\widetilde{\rA}^{(\nu,\kappa)}_{22}(s,x;t,y;T,t_1)$, 
$0 <s,t<T, (x,y) \in (0, \infty)^2$
are given such that
\begin{eqnarray}
&&\rho^{N \delta_0}_{T,(\nu,\kappa)} (t_{1}, \x^{(1)}_{N_1} ; 
\dots; t_M, \x^{(M)}_{N_M}) 
\nonumber\\
&&
=\Pf_{
\substack{1 \leq j \leq N_{m}, 1 \leq k \leq N_{n} \\ 1 \leq m, n \leq M}
}
\left[ \left(
\begin{array}{ll}
\widetilde{\rA}^{(\nu,\kappa)}_{11}(t_m,x^{(m)}_j; t_n, x^{(n)}_k;T,t_1) 
& \widetilde{\rA}^{(\nu,\kappa)}_{12}(t_m, x^{(m)}_j;t_n, x^{(n)}_k;T,t_1) 
\cr
-\widetilde{\rA}^{(\nu,\kappa)}_{12}(t_n, x^{(n)}_k;t_m, x^{(m)}_j;T,t_1) 
& \widetilde{\rA}^{(\nu,\kappa)}_{22}(t_m, x^{(m)}_j;t_n, x^{(n)}_k;T,t_1)
\end{array} \right) \right].
\nonumber
\end{eqnarray}
By direct calculation we have found that
\begin{eqnarray}
&& \widetilde{\rA}^{(\nu,\kappa)}_{jk}
(\sigma^2 c_{\sigma^2}(s), c_{\sigma^2}(s)^2x;
\sigma^2 c_{\sigma^2}(t),c_{\sigma^2}(t)^2y; \sigma^2, c_{\sigma^2}(t_M))
\nonumber\\
&& \quad
= c_{\sigma^2}(s)^{-1} c_{\sigma^2}(t)^{-1}
\rA^{(\nu,\kappa)}_{jk}(s,x;t,y;\sigma^2),
\label{eqn:AAnu}
\end{eqnarray}
for $(j,k)=(1,1), (1,2), (2,2)$,
where $\rA^{(\nu,\kappa)}_{jk}, (j,k)=(1,1), (1,2), (2,2)$, 
are given by (\ref{eqn:Anus}) with (\ref{eqn:BCnu}).
Then by the equality (\ref{eqn:Eq2}) of Lemma \ref{thm:equivalence}
and the property of Pfaffian (\ref{eqn:Pfaffian2}),
the theorem is proved. \qed

\vskip 0.5cm
\begin{small}
%%%%%%%%%%%%%%%%%%%%%%%%%%%%%%%%%%%%%%%%%%%%%%%%%%%%%%
\noindent{\bf Acknowledgements} \quad
%%%%%%%%%%%%%%%%%%%%%%%%%%%%%%%%%%%%%%%%%%%%%%%%%%%%%%
The present author would like to thank
T. Imamura for useful comments on the manuscript.
A part of the present work was done
during the participation of the present author 
in the ESI program ``Combinatorics and Statistical Physics"
(March and May in 2008).
The author expresses his gratitude for 
hospitality of the Erwin Schr\"odinger Institute 
(ESI) in Vienna
and for well-organization of the program
by M. Drmota and C. Krattenthaler.
This work is supported in part by
the Grant-in-Aid for Scientific Research (C)
(No.21540397) of Japan Society for
the Promotion of Science.
%%%%%%%%%%%%%%%%%%%%%%%%%%%%%%%%%%%%%%%%%%%%%%%%%%

%%%%%%%%%%%%%%%%%%%%%%%%%%%%%%%%%%%%%%%%%%%%%%%%%%%%%%%%%%%
%%%%%REFERENCES%%%%%%%%%%%%%%%%%%%%%%%%%%%%%%%%%%%%%%%%%%%%%
%%%%%%%%%%%%%%%%%%%%%%%%%%%%%%%%%%%%%%%%%%%%%%%%%%%%%%%%%%%%
%\begin{small}

%%%%%%%%%%%%%%%%%%%%%%%%%%%%%%%%%%%%%%%%%%%%%%%%%%%%%%%%%%
\end{small}
%%%%%%%%%%%%%%%%%%%%%%%%%%%%%%%%%%%%%%%%%%%%%%%%%%%%%%%%%%
\end{document}